\documentclass[12pt]{amsart}
\usepackage{amscd,amssymb,epsfig}
\usepackage[mathscr]{eucal}
\usepackage[english]{babel}
\usepackage{psfrag}

\input xy
\usepackage[all]{xy}

\begin{document}


\newtheoremstyle{mytheorem}
  {}
  {}
  {\slshape}
  {}
  {\scshape}
  {.}
  { }
  {}

\newtheoremstyle{mydefinition}
  {}
  {}
  {\upshape}
  {}
  {\scshape}
  {.}
  { }
  {}

\theoremstyle{mytheorem}
\newtheorem{lemma}{Lemma}[section]
\newtheorem{prop}[lemma]{Proposition}
\newtheorem{prop_intro}{Proposition}
\newtheorem{cor}[lemma]{Corollary}
\newtheorem{cor_intro}[prop_intro]{Corollary}
\newtheorem{thm}[lemma]{Theorem}
\newtheorem{thm_intro}[prop_intro]{Theorem}
\newtheorem*{thm*}{Theorem}
\theoremstyle{mydefinition}
\newtheorem{rem}[lemma]{Remark}
\newtheorem*{claim*}{Claim}
\newtheorem{rem_intro}[prop_intro]{Remark}
\newtheorem{rems_intro}[prop_intro]{Remarks}
\newtheorem*{notation*}{Notation}
\newtheorem*{warning*}{Warning}
\newtheorem{rems}[lemma]{Remarks}
\newtheorem{defi}[lemma]{Definition}
\newtheorem*{defi*}{Definition}
\newtheorem{defi_intro}[prop_intro]{Definition}
\newtheorem{defis}[lemma]{Definitions}
\newtheorem{exo}[lemma]{Example}
\newtheorem{exo_intro}[prop_intro]{Example}
\newtheorem{exos_intro}[prop_intro]{Examples}
\newtheorem*{exo*}{Example}
\newtheorem*{que*}{Question}

\numberwithin{equation}{section}

\newcommand{\bqn}{\begin{eqnarray*}}
\newcommand{\eqn}{\end{eqnarray*}}
\newcommand{\bq}{\begin{eqnarray}}
\newcommand{\eq}{\end{eqnarray}}
\newcommand{\ba}{\begin{aligned}}
\newcommand{\ea}{\end{aligned}}
\newcommand{\be}{\begin{enumerate}}
\newcommand{\ee}{\end{enumerate}}

\newcommand{\bibURL}[1]{{\unskip\nobreak\hfil\penalty50{\tt#1}}}

\def\ti{-\allowhyphens}

\newcommand{\thismonth}{\ifcase\month 
  \or January\or February\or March\or April\or May\or June%
  \or July\or August\or September\or October\or November%
  \or December\fi}
\newcommand{\thismonthyear}{{\thismonth} {\number\year}}
\newcommand{\thisdaymonthyear}{{\number\day} {\thismonth} {\number\year}}

\newcommand{\CC}{{\mathbb C}}
\newcommand{\DD}{{\mathbb D}}
\newcommand{\FF}{{\mathbb F}}
\newcommand{\HH}{{\mathbb H}}
\newcommand{\GG}{{\mathbb G}}
\newcommand{\KK}{{\mathbb K}}
\newcommand{\LL}{{\mathbb L}}
\newcommand{\NN}{{\mathbb N}}
\newcommand{\PP}{{\mathbb P}}
\newcommand{\QQ}{{\mathbb Q}}
\newcommand{\RR}{{\mathbb R}}
\renewcommand{\SS}{{\mathbb S}}
\newcommand{\TT}{{\mathbb T}}
\newcommand{\ZZ}{{\mathbb Z}}

\newcommand{\Bb}{{\mathcal B}}
\newcommand{\Cc}{{\mathcal C}}
\newcommand{\Dd}{{\mathcal D}}
\newcommand{\Ff}{{\mathcal F}}
\newcommand{\Ee}{{\mathcal E}}
\newcommand{\Gg}{{\mathcal G}}
\newcommand{\Hh}{{\mathcal H}}
\newcommand{\Kk}{{\mathcal K}}
\newcommand{\Ll}{{\mathcal L}}
\newcommand{\Mm}{{\mathcal M}}
\newcommand{\Nn}{{\mathcal N}}
\newcommand{\Oo}{{\mathcal O}}
\newcommand{\Pp}{{\mathcal P}}
\newcommand{\Qq}{{\mathcal Q}}
\newcommand{\Rr}{{\mathcal R}}
\newcommand{\Ss}{{\mathcal S}}
\newcommand{\Tt}{{\mathcal T}}
\newcommand{\Vv}{{\mathcal V}}
\newcommand{\Xx}{{\mathcal X}}
\newcommand{\Yy}{{\mathcal Y}}
\newcommand{\Zz}{{\mathcal Z}}

\newcommand{\fraka}{{\mathfrak a}}
\newcommand{\frakg}{{\mathfrak g}}
\newcommand{\frakk}{{\mathfrak k}}
\newcommand{\frakl}{{\mathfrak l}}
\newcommand{\frakp}{{\mathfrak p}}
\newcommand{\fraks}{{\mathfrak s}}
\newcommand{\fraku}{{\mathfrak u}}
\newcommand{\frakB}{{\mathfrak B}}

\newcommand{\dD}{{\mathbf D}}
\newcommand{\gG}{{\mathbf G}}
\newcommand{\hH}{{\mathbf H}}
\newcommand{\pP}{{\mathbf P}}
\newcommand{\lL}{{\mathbf L}}
\newcommand{\qQ}{{\mathbf Q}}

\newcommand{\G}{{\Gamma}}
\newcommand{\g}{{\gamma}}
\renewcommand{\l}{{\lambda}}
\newcommand{\T}{{\mathrm{T}}}
\newcommand{\SL}{{\mathrm{SL}}}
\newcommand{\PSL}{{\mathrm{PSL}}}
\newcommand{\GL}{{\mathrm{GL}}}
\newcommand{\Sp}{{\mathrm{Sp}}}
\newcommand{\SU}{{\mathrm{SU}}}
\newcommand{\SO}{{\mathrm{SO}}}
\newcommand{\Spin}{{\mathrm{Spin}}}
\newcommand{\U}{{\mathrm{U}}}
\newcommand{\is}{{\mathrm{Is}}}
\newcommand{\length}{{\operatorname{length}}}
\newcommand{\Mod}{{\mathbf{Mod}}}
\newcommand{\<}{\langle}
\renewcommand{\>}{\rangle}
\newcommand{\ol}{\overline}
\renewcommand{\phi}{{\varphi}}

\def\binfty{\mathcal B^\infty_{\mathrm {alt}}}
\def\bu{\bullet}
\def\cb{{\rm C}_{\rm b}}
\def\ehbc{{\rm EH}_{\rm cb}}
\def\Ef{\mathcal E_\varphi}
\def\essim{\operatorname{EssIm}}
\newcommand{\essimfi}{\operatorname{EssIm(\varphi)}}
\newcommand{\esssup}{\operatorname{ess\,sup}}
\def\h{{\rm H}}
\def\hb{{\rm H}_{\rm b}}
\def\hc{{\rm H}_{\rm c}}
\def\hcb{{\rm H}_{\rm cb}}
\def\hom{\operatorname{Hom}}
\def\homeo#1{{\sl H\!omeo}^+\!\left(#1\right)}
\def\thomeo#1{\widetilde{{\sl \!H}\!omeo}^+\!\left(#1\right)}
\def\id{{\it I\! d}}
\def\ind{\mathrm{ind}}
\def\la{\mathrm{L}^\infty_{\mathrm{alt}}}
\def\linfty{\mathrm{L}^\infty}
\def\linftyw{\mathrm{L}^\infty_{\rm w*}}
\def\lp{\mathrm{L}^p}
\def\ltwo{\mathrm{L}^2}
\def\rep{\operatorname{Rep}}
\newcommand{\supp}{\operatorname{supp}}

\def\one{\mathbf{1\kern-1.6mm 1}}
\def\property{\textbf{\rm\textbf A}}
\def\oddex#1#2{\left\{#1\right\}_{o}^{#2}}
\def\comp#1{{\rm C}^{(#1)}}
\def\lra{\longrightarrow}

\def\adg{\operatorname{ad}_\frakg}
\def\bg{B_\frakg}
\def\cs{\check S}
\def\deta{{\operatorname{det}_A}}
\def\gl{\mathrm{GL}}
\def\gg{\Gamma_g}
\def\ghtp{generalized Hermitian triple product }
\def\gmodp{\gG(\RR)/\pP(\RR)}
\def\gmodq{\gG(\RR)/\qQ(\RR)}
\def\gr{\mathrm{Gr}_p(W)}
\def\gra{\mathrm{Gr}^A_p(V)}
\def\Gr{\mathrm{Gr}}
\def\htp{Hermitian triple product }
\newcommand{\Iso}{\operatorname{Iso}}
\def\Isom{\operatorname{Isom}}
\def\isp{\operatorname{Is}_{\<\cdot,\cdot\>}}
\def\isptwo{\operatorname{Is}_{\<\cdot,\cdot\>}^{(2)}}\def\ispth{\operatorname{Is}_{\<\cdot,\cdot\>}^{(3)}}
\def\isf{\operatorname{Is}_F}
\def\isfi{\operatorname{Is}_{F_i}}
\def\isft{\operatorname{Is}_F^{(3)}}
\def\isfit{\operatorname{Is}_{F_i}^{(3)}}
\def\isftwo{\operatorname{Is}_F^{(2)}}
\def\lin{\operatorname{Lin}}
\def\ll{{\Ll_1,\Ll_2}}
\def\kahler{K\"ahler }
\def\kg{\kappa_G}
\def\kgb{\kappa_G^{\rm b}}
\def\kib{\kappa_i^{\rm b}}
\def\kibt{\tilde\kappa_i^{\rm b}}
\def\kx{\kappa_\Xx}
\def\kxb{\kappa_\Xx^\mathrm{b}}
\def\ox{\omega_\Xx}
\def\po{\mathrm{PO}}
\def\pr{\operatorname{pr}}
\def\psl{\mathrm{PSL}}
\def\pu{\mathrm{PU}}
\def\pupq{\mathrm{PU}(p,q)}
\def\puvi{\mathrm{PU}\big(V,\<\cdot,\cdot\>_i\big)}
\def\rkx{\operatorname{rk}_\Xx}
\def\rky{\operatorname{rk}_\Yy}
\def\sg{\Sigma_g}
\def\sltwo{\mathrm{SL}(2,\RR)}
\def\slv{\mathrm{SL}(V)}
\def\sp{\mathrm{Sp}}
\def\stab{\operatorname{Stab}}
\def\supq{\mathrm{SU}(p,q)}
\def\su{\mathrm{SU}}
\def\suq{\mathrm{SU}(q,1)}
\def\suw{\mathrm{SU}(W)}
\def\suvi{{\rm SU}\big(V,\<\cdot,\cdot\>_i\big)}
\def\tr{\mathrm{T}_\rho}
\def\u{\mathrm{U}}
\def\vol{\operatorname{vol}}
\def\bsl{\backslash}
\def\eps{\epsilon}

\def\hom{\operatorname{Hom}}
\def\sym{\operatorname{Sym}}
\def\en{\operatorname{End}}
\def\aut{\operatorname{Aut}}
\def\inn{\operatorname{Inn}}
\def\out{\operatorname{Out}}
\def\rep{\operatorname{Rep}}
\def\diff{\operatorname{Diff}}
\def\isom{\operatorname{Isom}}
\def\tr{\operatorname{tr}}
\def\d{\operatorname{d}}

\title[Action of mapping class group]{The action of the mapping class group on maximal representations}
\author[A.~Wienhard]{Anna Wienhard}
\email{wienhard@math.ias.edu}
\address{School of Mathematics, Institute for Advanced Study, 1
  Einstein Drive, Princeton NJ 08540, USA}
\address{Department of Mathematics, University of Chicago, 5734 University Avenue, 
Chicago, IL 60637-1514, 
USA} 
\thanks{The author was partially supported by the Schweizer Nationalfond under PP002-102765 and by the National Science Foundation under agreement No. DMS-0111298.}
\date{\today}

\begin{abstract}
Let $\G_g$ be the fundamental group of a closed oriented Riemann
surface $\Sigma_g$, $g\geq2$, and let $G$ be a simple Lie group of
Hermitian type. The Toledo invariant defines the 
subset of maximal representations $\rep_{max}(\G_g, G)$ in the representation 
variety $\rep(\G_g, G)$. $\rep_{max}(\G_g, G)$
is a union of connected components with similar properties as
Teichm\"uller space $\Tt(\Sigma_g) =  \rep_{max}(\G_g,
\PSL(2,\RR))$. 
We prove that the mapping class group $\Mod_{\Sigma_g}$ 
 acts properly on  
$\rep_{max}(\G_g, G)$ when $G=  \Sp(2n,\RR)$, $\SU(n,n)$, $\SO^*(4n)$, $\Spin(2,n)$.
\end{abstract}
\maketitle

\vskip2cm

\setcounter{tocdepth}{1}
\tableofcontents

\section{Introduction}

Let $\G_g$ be the fundamental group of a closed oriented surface $\Sigma_g$ 
of genus $g\geq 2$. Let $G$ be a connected semisimple Lie group 
and $\hom(\G_g, G)$ the space of homomorphisms $\rho: \G_g \to G$. 
The automorphism groups of $\G_g$ and $G$ act on $\hom(\G_g, G)$ by
\bqn
\aut(\G_g)\times \aut (G) \times \hom(\G_g ,G) &\to & \hom (\G_g, G)\\
(\psi, \alpha, \rho)&\mapsto& \alpha \circ\rho \circ\psi^{-1}:\left(\g \mapsto \alpha(\rho(\psi^{-1}\g))\right).
\eqn
Considering homomorphisms only up to conjugation in $G$ defines the representation variety 
\bqn
\rep(\G_g, G) := \hom(\G_g, G) / \inn(G).
\eqn
The above action induces an action of the group of outer automorphisms 
$\out(\G_g):= \aut(\G_g)/\inn(\G_g)$ of $\G_g$ 
on $\rep(\G_g, G)$: 
\bqn
\out(\G_g) \times \rep(\G_g, G) &\to& \rep(\G_g, G) \\
(\psi, [\rho]) &\mapsto& [\psi\rho] := \left[\left(\g \mapsto \rho(\psi^{-1} \g)\right)\right].
\eqn

Recall that $\out(\G_g)$ is isomorphic to $ \pi_0 (\diff(\Sigma_g))$. 
The mapping class group $\Mod_{\Sigma_g}$ is the subgroup of
$\out(\G_g)$ 
corresponding to orientation preserving diffeomorphisms of
$\Sigma_g$. 
We refer to \cite{Ivanov_mapping, Farb_Margalit} for a general introduction to mapping class groups 
and to \cite{Goldman_survey} for a recent survey on dynamical properties of 
the action of $\out(\G_g)$ on representation varieties $\rep(\G_g, G)$.

This note is concerned with the action of the mapping class group
on 
special connected components of $\rep(\G_g, G)$ 
when $G$ is of Hermitian type. 
Recall that a connected semisimple Lie group $G$ with finite 
center is said to be of Hermitian type if 
its associated symmetric space $\Xx$ is a Hermitian symmetric
space.
When $G$ is of Hermitian type there exists a bounded continuous integer valued function 
\bqn
\T: \rep(\G_g,G) \to \ZZ
\eqn
called the {\em Toledo invariant}. 
 
The level set of the maximal possible modulus of $\T$ is the set of 
{\em maximal representations}
\bqn
\rep_{max}(\G_g, G) \subset \rep(\G_g, G), 
\eqn
which is studied in \cite{Goldman_thesis, Goldman_88, Toledo_89,
  Hernandez, 
Bradlow_GarciaPrada_Gothen, Gothen,  
  Burger_Iozzi_Wienhard_ann, Burger_Iozzi_Labourie_Wienhard,
  Burger_Iozzi_Wienhard_tol, GarciaPrada_Gothen_Mundet}. 
Since the Toledo invariant is locally constant, its level sets are unions of connected components.

Results of \cite{Goldman_thesis, Goldman_88, Burger_Iozzi_Wienhard_ann,
  Burger_Iozzi_Wienhard_tol}  suggest that maximal representations
  provide a meaningful generalization of 
Teichm\"uller space when  $G$ is of Hermitian type \cite{Wienhard_grenoble}.
This note supports this similarity by proving the following theorem
\begin{thm}\label{thm:intro_main}
Let $G= \Sp(2n,\RR), \, \SU(n,n),\, \SO^*(4n),\,\Spin(2,n)$. 
Then the action of $\Mod_{\Sigma_g}$ on $\rep_{max}(\G_g, G)$ is proper.
\end{thm}
The validity of Theorem~\ref{thm:intro_main} for all groups locally
isomorphic to either 
$\Sp(2n,\RR)$, $\SU(n,n)$, $\SO^*(4n)$ or $\Spin(2,n)$ 
would follow from an affirmative answer to the following question: 
\begin{que*}
If $G= \Sp(2n,\RR)$, $\SU(n,n)$, $\SO^*(4n)$ or $\Spin(2,n)$, $\ol{G}$
the adjoint form of $G$, and $\rho\in \rep_{max}(\G_g, \ol{G})$, does
there exist a lift of $\rho$ to G?
\end{que*}

\begin{rem}
Note that maximal representations factor through maximal subgroups of
tube type \cite{Burger_Iozzi_Wienhard_ann, Burger_Iozzi_Wienhard_tol}. 
Therefore the only case which is not covered by the above theorem is the exceptional group $G={E_7}_{(-25)}$.
\end{rem}

We would like to remark that 
the study of maximal representations $\rep_{max}(\G_g, G) \subset \rep(\G_g, G)$ when $G$ is of Hermitian type is related to 
the study of
the Hitchin component $\rep_H(\G_g, G) \subset \rep(\G_g, G)$ for split real
simple Lie groups $G$. Fran\c{c}ois Labourie recently announced, as a
consequence of his work on Anosov representations and crossratios 
\cite{Labourie_anosov, Labourie_crossratio}, that the
mapping class group acts properly on $\rep_H(\G_g, \SL(n,\RR))$. After finishing this note, we learned that he also has a proof for maximal representations into $\Sp(2n,\RR)$ \cite{Labourie_energy}.

The author is indebted to Marc Burger for motivation, interesting  discussions and for pointing out 
a mistake in a preliminary version of this paper. The author thanks 
Bill Goldman, Ursula Hamenst\"adt, Alessandra Iozzi and Fran\c{c}ois
Labourie for useful discussions, and the referee for detailed
suggestions which helped to substantially improve the exposition of this note. 

\vskip1cm
\section{Maximal Representations and Translation Lengths}
\subsection{Maximal Representations}
For an introduction and overview the reader is referred to
\cite{Burger_Iozzi_Labourie_Wienhard, Burger_Iozzi_Wienhard_tol}. 
Let $G$ be a connected semisimple Lie group with finite center. Denote
by $\Xx= G/K$, with $K< G$ a maximal compact subgroup, its
associated symmetric space.  
$G$ is said to be of {\em Hermitian type} if there exists a $G$-invariant complex structure 
$J$ on $\Xx$. The composition of the Riemannian metric induced by the Killing form $\frakB$ on
$\Xx$ with the complex structure 
defines a K\"ahler form 
\bqn
\omega_\Xx(v,w):= \frac{1}{2}\frakB(v, Jw) 
\eqn
which is a $G$-invariant closed differential
two-form on $\Xx$. 

Given a representation $\rho: \G_g \to G$ consider the 
associated flat bundle $E_\rho$ over $\Sigma_g$ 
defined by 
\bqn
E_\rho:= \G_g\bsl(\widetilde\Sigma_g\times \Xx), 
\eqn
where $\G_g$ acts diagonally by deck transformations on
$\widetilde\Sigma_g$ and via $\rho$ on $\Xx$.  
As $\Xx$ is contractible, there
exists a smooth section  $f: \Sigma_g\to E_\rho$ which is unique up to
homotopy. 
This section lifts to a smooth $\rho$-equivariant map $\tilde{f}:
\widetilde\Sigma_g \to \widetilde\Sigma_g \times \Xx \to \Xx$.
The pull back of $\omega_\Xx$ via 
$\tilde{f}$ is a $\G_g$-invariant two-form $\tilde{f}^* \omega_\Xx$ on
$\widetilde\Sigma_g$ which may 
be viewed as a two-form on the closed
surface $\Sigma_g$. 
The {\em Toledo invariant} of $\rho$ is  
\bqn
\T(\rho):= \frac{1}{2\pi} \int_{\Sigma_g} \tilde{f}^* \omega_\Xx.
\eqn
The Toledo invariant is independent of the choice of the section $f$
and defines a continuous function
\bqn
\T:\hom(\G_g, G) \to \ZZ.
\eqn
The map $\T$ is invariant under the action of $\inn(G)$ and constant on connected components
of the representation variety. 
The Toledo invariant satisfies a generalized Milnor-Wood inequality \cite{Domic_Toledo, Clerc_Orsted_2}
\bqn
|\T|\leq\frac{p_\Xx\rkx}{2}| \chi(\Sigma_g)|,
\eqn
where $\rkx$ is the real rank of $\Xx$ and $p_\Xx \in \NN$ 
is explicitly computable in terms of the root system.

\begin{defi}
A representation $\rho:\G_g\to G$ is said to be {\em maximal} if 
\bqn
|\T(\rho)| = \frac{p_\Xx \rkx}{2}|\chi(\Sigma_g)|.
\eqn
\end{defi}
\begin{rem}
Changing the orientation of $\Sigma_g$ switches the sign of $\T$. We will restrict our attention to the 
case when $\rho$ is maximal with $\T(\rho) >0$. 
\end{rem}

We define the set of maximal representations 
\bqn
\rep_{max}(\G_g, G) := \{[\rho]\in \rep(\G_g, G) \, |\, \rho \text{ is
  a maximal representation}\}, 
\eqn
which is a union of connected components of $
\rep(\G_g, G)$.
The set $\rep_{max}(\G_g, \PSL(2,\RR))$ is the union of the two
Teichm\"uller components of $\Sigma_g$ \cite{Goldman_thesis}.

The action of the group $\out(\G_g):= \aut(\G_g)/\inn(\G_g)$  of outer
automorphism of $\G_g$ 
on $\rep(\G_g, G)$ given by 
\bqn
\out(\G_g) \times \rep(\G_g, G) &\to& \rep(\G_g, G) \\
(\psi, [\rho]) &\mapsto& [\psi\rho] := [(\g \mapsto \rho(\psi \g))].
\eqn
preserves $\rep_{max}(\G_g, G)$. 

The mapping class group $\Mod_{\Sigma_g}$ preserves, and hence acts on
the components of $\rep_{max}(\G_g, G)$ where $\T>0$.
\begin{rem}
Note that whereas Teichm\"uller space, the set of quasifuchsian
representations and Hitchin components are always contractible
subsets of $\rep(\G_g, G)$, certain
components of the set of maximal representations might have nontrivial
topology \cite{Gothen, Bradlow_GarciaPrada_Gothen_homotopy}.
\end{rem}
%
\subsection{Translation Lengths}
For a hyperbolization $h:\G_g\to \PSL(2,\RR)$ define the {\em  translation length} $\tr_h$ of $\g\in \Gamma_g$ as  
\bqn
\tr_h(\g):= \inf_{p\in \DD} \d_\DD(p, \g p).
\eqn
For a representation $\rho: \Gamma_g\to G$ define similarly the
  translation length $\tr_\rho$ of $\g\in \G_g$ as 
\bqn
\tr_\rho(\g):= \inf_{x\in \Xx_G} \d_\Xx (x, \rho(\g) x), 
\eqn
where $\d_\Xx$ is any left-invariant distance on 
the symmetric space associated to $G$. 

\begin{prop}\label{prop:trlength}
Fix a hyperbolization $h$ of $\G_g$. Assume that for any maximal
representation $\rho: \G_g\to G$ there exists $A,B>0$ such that  
\bq\label{eq:compatible} 
A^{-1} \tr_h(\g) - B \leq \tr_{\rho}(\g) \leq A \tr_h(\g) +B \quad \text{ for all $\g\in \G_g$.}
\eq
Then $\Mod_{\Sigma_g}$ acts properly on $\rep_{max}(\G_g, G)$.
\end{prop}
The Proposition relies on the fact that $\Mod_{\Sigma_g}$ acts properly
discontinuous on Teichm\"uller space $\Tt(\G_g)$, which is due to
Fricke. 
\begin{lemma}\cite[Proposition 5]{Douady}\label{lem:finite_curves}
There exists a collection of simple closed curves $\{c_1, \cdots c_{9g-9}\}$ on $\Sigma_g$ 
such that the map 
\bqn
\Tt(\G_g) &\to& \RR^{9g-9}\\
h &\mapsto& (\tr_h(\g_i))_{i=1,\cdots, 9g-9}, 
\eqn
where $\g_i$ is the element of $\G_g$ corresponding to $c_i$, 
is injective and proper. 
\end{lemma}
\begin{rem}
A family of such $9g-9$ curves is given by $3g-3$ curves $\alpha_i$ giving a pants decomposition,
$3g-3$ curves $\beta_i$ representing seems of the pants decomposition and 
the $3g-3$ curves given by the Dehn twists of $\beta_i$ along $\alpha_i$ (see e.g. \cite{Farb_Margalit}). 
\end{rem}

\begin{proof}[Proof of Proposition~\ref{prop:trlength}]
We argue by contradiction. Suppose that the action of
$\Mod_{\Sigma_g}$ on $\rep_{max}(\G_g, G)$ is
not proper. Then there exists a compact subset
$C\subset \rep_{max}(\G_g, G)$ such that 
\bqn
\#\{\psi \in \Mod_{\Sigma_g} \, |\,
\psi(C) \cap C\} 
\eqn
is infinite. Thus there exists an infinite sequence $\psi_n$ 
in $\Mod_{\Sigma_g}$ and a representation $\rho\in \rep_{max}(\G_g, G)$ such that 
$\psi_n(\rho)$ converges to a representation $\rho_\infty\in
\rep_{max}(\G_g, G)$. 
Since $\psi_n$ acts properly on Teichm\"uller space $\Tt(\G_g)$, 
the sequence
of hyperbolizations $\psi_n h$ leaves every compact set of
$\Tt(\G_g)$. 
This implies that the sum of the translation lengths of the 
elements $\g_i$, $ i= 1, \cdots, 9g-9$ tends to $\infty$: 
\bqn
\sum_{i=1}^{9g-9} \tr_{\psi_n h} (\g_i) \to \infty
\eqn
By assumption (\ref{eq:compatible}) 
\bqn
A^{-1} \tr_h(\psi_n^{-1}\g_i) - B \leq \tr_{\rho}(\psi_n^{-1}\g_i), 
\eqn   
hence 
\bqn
\sum_{i=1}^{9g-9} \tr_{\psi_n \rho} (\g_i) \to \infty.
\eqn
This contradicts $\lim_{n\to \infty}\psi_n\rho = \rho_\infty$, since, by (\ref{eq:compatible}), the sum 
$\sum_{i=1}^{9g-9} \tr_{\rho_\infty} (\g_i)$ is bounded from above by 
$A \sum_{i=1}^{9g-9} \tr_{h} (\g_i)+ B$. 
\end{proof}
Note that the upper bound for the comparison of the translation lengths with respect to a hyperbolization $h$ and 
a representation $\rho$
is established quite easily
\begin{lemma}\label{lem:upperbound}
Fix a hyperbolization $h$. 
For every maximal representation $\rho: \G_g \to G$ there exists $A,B\geq 0$
such that  
\bqn
\tr_{\rho}(\g) \leq A\tr_h(\g) + B \quad \text{ for all $\g \in \G_g$}.
\eqn
\end{lemma}
\begin{proof}
Let $\Xx$ be the symmetric space associated to $G$.
By \cite[Proposition~2.6.1]{Korevaar_Schoen} there exists a $\rho$-equivariant 
(uniform) $L$-Lipschitz map $f: \DD \to \Xx$.
Let $p_0\in \DD$ such that $\tr_h(\g) = \d_\DD(p_0, \g p_0)$, then
\bqn
\tr_\rho(\g) &\leq& \d_\Xx(f(p_0), \rho(\g) f(p_0)) = \d_\Xx(f(p_0), f(\g p_0))\\
&\leq& L \d_\DD(p_0, \g p_0) = L \tr_h(\g).
\eqn
\end{proof}

\vskip1cm
\section{Maximal Representations into the Symplectic Group}
The main objective of this section is to establish the following

\begin{prop}\label{prop:symp}
For any hyperbolization $h$ of $\G_g$, there exist constants
$A,B\geq 0$ such that  
\bqn
A^{-1} \tr_h(\g) - B \leq \tr_{\rho}(\g) \leq A \tr_h(\g) +B 
\eqn
for all $\rho \in \rep_{max}\left(\G_g, \Sp(2n,\RR)\right)$ and all $\g\in \G_g$.
\end{prop}
Proposition~\ref{prop:symp} in combination with
Proposition~\ref{prop:trlength} gives 
\begin{cor}
The action of $\Mod_{\Sigma_g}$ on $\rep_{max}(\G_g, \Sp(2n,\RR))$ is proper. 
\end{cor}

That Theorem~\ref{thm:intro_main} can be deduced from
Proposition~\ref{prop:symp} and Proposition~\ref{prop:trlength} can be
seen as follows - we refer the reader to
\cite{Burger_Iozzi_Labourie_Wienhard, Burger_Iozzi_Wienhard_tight,
  Wienhard_thesis} for more on tight homomorphisms and their
properties. 
Satake \cite[Ch.~IV]{Satake_book} investigated when a simple Lie
group $G$ of Hermitian type admits a homomorphism 
\bqn
\tau: G\to \Sp(2m, \RR).
\eqn
such that the induced homomorphism of Lie algebras 
\bqn
\pi: \frakg\to \fraks\frakp(2m,\RR)
\eqn
is a so called 
 $(H_2)$-Lie algebra homomorphism. Examples of such are   
\bqn
&\tau:& \SU(n,n) \to \Sp(4n, \RR)\\
&\tau:& \SO^*(4n) \to \Sp(8n,\RR)\\
&\tau:& \Spin(2,n) \to \Sp(2m, \RR), \text{ where $m$ depends on }n\mod 8.
\eqn

In \cite{Wienhard_thesis, Burger_Iozzi_Wienhard_tight} we prove that any
such $(H_2)$-homomorphism $\tau$ is a tight homomorphism. This implies in particular that the
composition of any maximal representation $\rho: \G_g\to G$ for $G= \SU(n,n),\, \SO^*(4n),\,
\Spin(2,n)$ with the homomorphism $\tau: G\to \Sp(2m, \RR)$ is a
maximal representation $\rho_\tau:=\tau\circ \rho: \G_g \to \Sp(2m, \RR)$. 
By Proposition~\ref{prop:symp} the translation lengths $\tr_h(\g)$
and $\tr_{\rho_\tau}(\g)$ are comparable. 
Since the embedding $\Xx_G \to \Xx_{\Sp(2m,\RR)}$, defined by $\tau$, is totally geodesic 
and the image $\rho_\tau(\G_g)$ preserves $\Xx_G$, the same
argument as in Lemma~\ref{lem:subspace} below gives that $\tr_{\rho_\tau}(\g) = \tr_{\rho}(\g)$ for all
$\g\in {\G_g}$.
%
%
\subsection{The Symplectic Group}
For a $2n$-dimensional real vector space $V$ with a nondegenerate skew-symmetric
bilinear form $\langle \cdot, \cdot \rangle$, the 
symplectic group $\Sp(V)$ is defined as 
\bqn
\Sp(V):=\aut(V,\langle \cdot, \cdot \rangle) :=\{g \in \GL(V)\,
|\,\langle g\cdot, g\cdot \rangle=\langle  \cdot, \cdot \rangle \}.
\eqn

The symmetric space associated to $\Sp(V)$ is given by 
\bqn
\Xx_{\Sp}:= \{ J\in \GL(V) \, |\, J^2 = -\id, \, \,\langle  \cdot, J \cdot
\rangle >> 0 \}, 
\eqn
where $\langle  \cdot, J \cdot
\rangle >> 0$ indicates that $\langle  \cdot, J \cdot
\rangle$ is symmetric and positive definite. The action of $\Sp(V)$ on $\Xx_\Sp$ is 
by conjugation $g(J) = g^{-1}Jg$.

We specify a left invariant distance on $\Xx_\Sp$ as follows. 
Let $J_1, J_2 \in \Xx_\Sp$, the symmetric positive definite forms 
$\<\cdot, J_i \cdot\>$ define a pair of Euclidean norms $q_i$ on $V$. 
Denoting by $||\id||_{J_1, J_2}$ the norm of the identity map from
$(V, q_1)$ to $(V, q_2)$ we define a distance on $\Xx_\Sp$ by 
\bqn
\d_\Sp (J_1, J_2):= \Big{|}\log ||\id||_{J_1, J_2}\Big{|} +  \Big{|}\log ||\id||_{J_2, J_1}\Big{|} 
\eqn
 
\subsection{Transverse Lagrangians and Causal Structure}
Let 
\bqn
\Ll(V): = \{L\subset V\, |\, \dim(L)=n, \, \,{\langle  \cdot, \cdot
\rangle_|}_L = 0 \}
\eqn
be the space of Lagrangian subspaces of $V$. 
Two Lagrangian subspaces $L_+, L_- \in \Ll(V)$ are said to be
transverse if $L_+\cap L_- = \{0\}$. 
Any two transverse Lagrangian subspaces $L_+, L_-\in \Ll(V)$ define a
symmetric subspace $\Yy_{L_-, L_+} \subset \Xx_\Sp$ by 
\bqn
\Yy_{L_-, L_+}:= \{ J \in \Xx_\Sp\, |\, J(L_\pm) = L_\mp\} \subset \Xx_\Sp.
\eqn
Writing an element $g\in \Sp(V)$ in block decomposition 
$g=\begin{pmatrix} A& B\\C&D
\end{pmatrix}$ with
respect to the decomposition $V= L_-\oplus L_+$ defines a 
natural embedding $\GL(L_-)\to
\Sp(V)$ given by 
\bqn
\GL(L_-) &\to& \Sp(V)\\
A &\mapsto& \begin{pmatrix} A& 0\\0&{A^T}^{-1}
\end{pmatrix}, 
\eqn 
similarly for $\GL(L_+)$.

The subgroup $\GL(L_-)$ preserves the symmetric subspace $\Yy_{L_-, L_+}$
and 
acts transitively
on it. 

\begin{rem}
The space of Lagrangians $\Ll(V)$ 
can be identified with the Shilov boundary $\cs_\Sp$ of $\Xx_\Sp$, and realized inside the
visual boundary as the $G$-orbit of a specific maximal
singular direction. 
Two Lagrangians are transverse if and only if
the two corresponding points in the visual boundary can be joined by a
maximal singular geodesic $\g_{L_\pm}$.
The symmetric subspace $\Yy_{L_-, L_+}$ is the parallel set of $\g_{L_\pm}$, i.e. the set of points 
on flats containing
the geodesic $\g_{L_\pm}$; it is the noncompact symmetric space dual to $\Ll(V)\simeq \U(L_-) / {\rm O}(L_-)$. 
\end{rem}

For $J\in \Yy_{L_-, L_+}$  the restriction of $\<
\cdot, J\cdot \>$ to $L_-$ is a 
positive definite symmetric bilinear
form on $L_-$, and conversely, fixing $L_+$,  
any positive definite symmetric bilinear
form $Z$ on $L_-$ defines a complex structure $J\in \Yy_{L_-, L_+}$.
Therefore, the space  
\bqn
\Yy_{L_-,s}:=\{Z \, |\, Z \text{ positive bilinear form on } L_-\}.
\eqn
with the action of $\GL(L_-)$ by 
\bqn
A(Z):=
A^T Z{A}, 
\eqn 
where we choose a scalar product on $L_-$ (i.e. a base point) and realize a bilinear form on $L_-$ 
as a symmetric $(n\times n)$ matrix $Z$, 
is $\GL(L_-)$-equivariantly isomorphic to $\Yy_{L_-, L_+}$. 
We endow $\Yy_{L_-,s}$ with the left-invariant distance induced by $\d_\Sp$ via this isomorphism. 

The space $\Yy_{L_-, s}$ is endowed with a natural causal
structure, given by the $\GL(L_-)$-invariant family of proper open cones 
\bqn
\Omega_Z:= \{ Z' \subset \Yy_{L_-, s}\, |\, Z'-Z \text{ is positive definite }\}.
\eqn
\begin{defi}\label{defi:causal_cont}
A continuous map $f:[0,1] \to \Yy_{L_-, s}$ is said to be
\emph{causal} if $f(t_2) 
\in \Omega_{f(t_1)}$ for all $0\leq t_1< t_2\leq 1$.
\end{defi}
A consequence of the proof of Lemma~8.10 in \cite{Burger_Iozzi_Labourie_Wienhard} is the following 
\begin{lemma}\label{lem:causal_curve}
For all $Z\in \Yy_{L_-, s}$, $Z'\in \Omega_Z$  and every causal
  curve $f:[0,1] \to \Yy_{L_-, s}$ with  $f(0)=Z$ and $f(1)=Z'$: 
\bqn
\length(f) \leq n \d_\Sp(Z,Z'),
\eqn
where $n=\dim(L_-)$.
\end{lemma}
The claim basically follows from the last inequality in the proof of Lemma~8.10 in 
\cite{Burger_Iozzi_Labourie_Wienhard}. However, for the reader's
convenience we give a direct proof here. 
\begin{proof}
Since $\d_\Sp$ is left invariant it is enough to prove the statement for
$Z=\id_{n}\in \Yy_{L_-,s}$.
For any subdivision 
\bqn
0 = t_0 < t_1< \cdots <t_m = 1
\eqn
let $f(t_i) = B^T_i B_i \in \Yy_{L_-,s}$, and note that by causality  
\bqn
\det\left( (B_{i}B_{i+1}^{-1})^T (B_{i}B_{i+1}^{-1})\right)<1.
\eqn
With $n= \dim(L_-)$, we have
\bqn
\d_\Sp(f(t_i), f(t_{i+1})) &=& \d_\Sp(B_i^T B_i, B_{i+1}^T B_{i+1})\\
& =& \log\left[ \lambda_{max}\left( (B_{i+1}B_i^{-1})^T (B_{i+1}B_i^{-1})\right)\right]\\ 
&+& \log \left[\lambda_{min}\left( (B_{i+1}B_i^{-1})^T (B_{i+1}B_i^{-1})\right)\right]\\
&\leq& \log \left[\det\left( (B_{i+1}B_i^{-1})^T (B_{i+1}B_i^{-1})\right)\right] \\
&-& n \log\left[ \det\left( (B_{i}B_{i+1}^{-1})^T (B_{i}B_{i+1}^{-1})\right)\right]\\
&\leq& n \log \left[\lambda_{max}( B_{i+1}^T B_{i+1})\right]\\
& -& n\log\left[\lambda_{min}( B_{i}^T B_{i})\right] \\
&+& n\log\left[\lambda_{min}( B_{i+1}^T B_{i+1})\right]\\
&-& n\log\left[\lambda_{max}( B_{i}^T B_{i})\right]. 
\eqn
Summing over the subdivision we obtain 
\bqn
\length(f)&\leq&\sum_{i+1}^m \d_\Sp (f(t_i), f(t_{i+1})) \\
&\leq& n \left[\log\lambda_{max}(f(1)) + \log \lambda_{min}(f(1))\right]\\
& -&  n \left[\log\lambda_{max}(f(0)) -  \log \lambda_{min}(f(0))\right]\\
&=&  n \left[\log\lambda_{max}(f(1)) +  \log \lambda_{min}(f(1))\right] \\
&=& n \d_\Sp(f(0), f(1))
\eqn
as needed.
\end{proof}

\subsection{Quasi-isometric Embedding}
Let $\rho:\G_g \to \Sp(V)$ be a maximal representation. The choice of
a hyperbolization $h$ of $\G_g$ defines a natural action of $\G_g$ on
$S^1= \partial \DD$. 
\begin{lemma}\label{lem:map}\cite[Corollary~6.3]{Burger_Iozzi_Labourie_Wienhard}
There exists 
a $\rho$-equivariant continuous map 
$\phi: S^1\to \Ll(V)$ such that distinct points 
$x, y \in S^1$ are mapped to transverse Lagrangians $\phi(x), \phi(y) \in \Ll(V)$.
\end{lemma}

A triple $(L_-, L_0, L_+)\in\Ll(V)^3$ of pairwise transverse Lagrangians gives rise to a complex structure 
\bqn
J_{L_0}=\begin{pmatrix} 0 & -T^+_0\\ T^-_0 &0\end{pmatrix} \quad \text{on}\quad V=L_- \oplus L_+,
\eqn 
where $T^\pm_0:
  L_\pm \to L_\mp$ is the unique linear map such that $L_0 = {\operatorname{graph}}(T^\pm_0)$. 
A triple $(L_-, L_0, L_+)$ of pairwise transverse Lagrangians is {\em maximal} if the 
symmetric bilinear form $\<\cdot, J_{L_0}\cdot\>$ is positive definite, that is if 
$J_{L_0} \in \Yy_{L_-, L_+} \subset \Xx_\Sp$.
We denote by $\Ll(V)^{3_+}$ the space of maximal pairwise transverse triples in $\Ll(V)$. 

Under the identification of the unit tangent bundle of the
  Poincar\'e-disc $T^1\DD\simeq (S^1)^{3_+}$ with positively oriented
  triples in $S^1$, 
the map $\phi$ gives rise to a $\rho$-equivariant map (Equation $(8.9)$ in \cite{Burger_Iozzi_Labourie_Wienhard})
\bqn
J: T^1\DD\cong (S^1)^{3_+} \to \Ll(V)^{3_+} &\to& \Xx_\Sp\\
u= (u_-, u_0, u_+) \mapsto \left(\phi(u_-), \phi(u_0), \phi(u_+)\right) &\mapsto& J(u),
\eqn
where $J(u)$ is the complex structure defined by the maximal triple 
\bqn
\left(\phi(u_-), \phi(u_0), \phi(u_+)\right)\in \Ll(V)^{3_+}. 
\eqn
Let $g_t$ be the lift of the geodesic flow on $T^1\Sigma_g$ to $T^1\DD$. 
Then for all $t$ the image of $g_t u = (u_-, u_t, u_+)$ under $J$ is contained in the symmetric subspace 
$\Yy_{\phi(u_-), \phi(u_+)} \subset \Xx_\Sp$ associated to the two 
transverse Lagrangians $\phi(u_-), \phi(u_+)$.

\begin{lemma}\cite[Equation ~8.8]{Burger_Iozzi_Labourie_Wienhard}\label{lem:strongquasi}
Let $\g\in \G_g\bsl\{\id\}$ and $p\in \DD$ a point. 
Denote by $u\in T^1\DD$ the unit tangent vector at $p$ of the geodesic connecting $p$ and $\g p$. 
Then there exists a constant $A'>0$ such that 
\bqn
A'^{-1} \d_\DD(\g p, p) \leq \d_\Sp (J(u), \rho(\g) J(u)).
\eqn
\end{lemma}

\begin{rem}
Note that the statement of Lemma~\ref{lem:strongquasi} implies that
the action of $\Mod_{\Sigma_g}$ on the connected components of maximal Toledo invariant in $\hom(\G_g, G)$ 
is proper, but is not sufficient to deduce
Theorem~\ref{thm:intro_main}. 
The inequality
in Lemma~\ref{lem:strongquasi} is with respect to specific points 
in $\Xx_\Sp$, but to compare the translation lengths we have
to take infima on both sides of the inequality. There is
in general no direct way to compare the translation length of
$\rho(\g)$ with the displacement length of $\rho(\g)$ with respect to
a specific point $x\in \Xx_\Sp$. In our situation we will make use of the causal structure on $\Yy_{L_-, s}$ 
to compare the translation length of $\rho(\g)$ with the displacement length 
$\d_\Sp (J(u), \rho(\g) J(u))$.
\end{rem}

\subsection{Translation Length and Displacement Length} 
Fix a hyperbolization $h: \G_g\to \PSL(2,\RR)$. 
Let $\g \in \G_g\bsl \{\id\}$. Denote by $\g^+, \g^- \in S^1$ the
attracting, respectively repelling, fixed point of $\g$ and 
by $L_\pm \in \Ll(V)$ the images of $\g^\pm$ under the
$\rho$-equivariant boundary map $\phi: S^1\to \Ll(V)$. 
Let $\Yy_{L_-, L_+}\subset \Xx_\Sp$ be the symmetric subspace associated
to $L_+, L_-$. 
\begin{lemma}\label{lem:subspace}
The translation length of $\rho(\g)$ is attained on $\Yy_{L_-, L_+}$:
\bqn
\tr_\rho(\g) = \inf_{J'\in \Yy_{L_-, L_+}} \d_\Sp(J', \rho(\g) J').
\eqn
\end{lemma} 
\begin{proof}
The symmetric subspace $\Yy_{L_-, L_+}$ is a totally geodesic submanifold of $\Xx_\Sp$. Denote by 
$\pr_\Yy: \Xx_\Sp \to \Yy_{L_-, L_+}$ the nearest point projection onto $\Yy_{L_-, L_+}$. The projection $\pr_\Yy$ is distance decreasing.
Since the element $\rho(\g)$ stabilizes $L_\pm$, the projection
$\pr_\Yy$ onto $\Yy_{L_-, L_+}$ 
is $\rho(\g)$-equivariant. 
Thus, for every $J\in \Xx_\Sp$
\bqn
\d_\Sp\left(J, \rho(\g) J\right) \geq 
\d_\Sp\left(\pr_\Yy(J), \pr_\Yy(\rho(\g)J)\right) =
\d_\Sp\left(\pr_\Yy(J), \rho(\g)\pr_\Yy(J)\right).
\eqn
In particular 
\bqn
\inf_{J\in \Xx_\Sp} \d_\Sp(J, \rho(g) J) = \inf_{J' \in \Yy_{L_-, L_+}} \d_\Sp(J', \rho(g) J').
\eqn
\end{proof}

\begin{lemma}\label{lem:properties}
Let $p_0\in \DD$ be some point lying on the geodesic $c_\g$ connecting $\g^-$ to $\g^+$. 
Let $u= (\g^-, u_0, \g^+) \in T^1\DD$ be the unit vector tangent to $c_\g$ at $p_0$. 
Then
\begin{enumerate}
\item $\rho(\g) J(u) \in \Omega_{J(u)}$.
\item For any  $Z\in Y_{L_-,s}$ there exists an $N$ such that for all $m\geq N$: 
$\rho(\g)^m J(u)\in \Omega_Z$. 
\end{enumerate}
\end{lemma}
\begin{proof}
(1) The map $u \mapsto J(u)$ is $\rho$-equivariant, thus $\rho(\g) J(u)= J(v) \in \Yy_{L_-, s}$,  
with $v= (\g^-, \g u_0, \g^+)$.
Since $\g^+$ is the attracting fixed point of $\g$, the triple 
$(u_0, \g u_0, \g^+)$ is positively oriented, and $(\phi(u_0), \phi(\g u_0), \phi(\g^+))$ is a maximal triple. Then, by 
\cite[Lemma~8.2]{Burger_Iozzi_Labourie_Wienhard}, $\rho(\g) J(u)-J(u)= J(v) -J(u)$ 
is positive definite, 
thus $\rho(\g) J(u) \in \Omega_{J(u)}$.\\
(2) Let $\mu$ be the maximal eigenvalue of $Z\in \Yy_{L_-, s}$. 
It suffices to show that there exists some 
$N$ such that the minimal eigenvalue of $\rho(\g) ^N J(u)\in \Yy_{L_-,
  s}$ is bigger than $\mu$. 
Then $\rho(\g)^N J(u) \in \Omega_Z$ and, with statement (1), we have that  
$\rho(\g)^m J(u) \in \Omega_{\rho(\g)^{N} J(u)} \subset \Omega_Z$ 
for all $m>N$.
Note that $\g^m u_0 \to \g^+$ as $m\to \infty$ and, since $\phi$ is continuous,
$\rho(\g)^m \phi(u_0) \to \phi(\g^+)$.
Moreover $\rho(\g)^{i+1}J(u) - \rho(\g)^i J(u)$ is positive definite for all $i$.  
Hence, the eigenvalues of $\rho(\g)^i J(u)$ grow monotonically towards $\infty$. 
In particular, there exists $N$ such that the minimal eigenvalue of 
$\rho(\g) ^N J(u)$ is bigger than $\mu$.
\end{proof}
Combining Lemma~\ref{lem:properties} with Lemma~\ref{lem:causal_curve} we obtain 
\begin{lemma}\label{lem:lengths}
There exist constants $A'', B'' >0$ such that for every $Z\in \Yy_{L_-, s}$   
\bqn
\d_\Sp(J(u), \rho(\g) J(u)) \leq A'' \d_\Sp(Z, \rho(\g)Z) + B''.
\eqn
\end{lemma}
\begin{proof}
Fix $Z\in \Yy_{L_-, s}$. Choose, by Lemma~\ref{lem:properties}, $N$ 
big enough such that $Z':= \rho(\g)^N J(u)\in \Omega_Z$.  
By Lemma~\ref{lem:properties} there are causal, distance realizing curves $f_Z$ from $Z$ to $Z'$ and 
$f_i$ from $\rho(\g)^i
Z'$ to $\rho(\g)^{i+1} Z'$ for all $0\leq i$. 
For every $k\geq0$ the concatenation $f= f_{k-1} *\cdots *f_0*f_Z$ is a causal curve from 
$Z$ to $\rho(\g)^k Z' = \rho(\g)^{N+k} J(u)$.
Thus applying Lemma~\ref{lem:causal_curve} we get that for every $k\geq 0$  
\bqn
&&\d_\Sp(Z,Z') + k \d_\Sp(Z', \rho(\g) Z') \\
&=& \d_\Sp(Z,Z') +
\sum_{i=0}^{k-1}\d_\Sp\left(\rho(\g)^i Z', \rho(\g)^{i+1} Z'\right)\\
&=&\length (f) \\
&\leq& n\d_\Sp(Z, \rho(\g)^k Z') \\
&\leq& n\left[\d_\Sp(Z, \rho(\g)^k Z) + \d_\Sp (\rho(\g)^k Z, \rho(\g)^k Z') \right]\\
&\leq& n\left[\sum_{i=0}^{k-1}\d_\Sp(\rho(\g)^i Z, \rho(\g)^{i+1} Z) + \d_\Sp (Z', Z) \right]\\
&=&  n \d_\Sp(Z',Z)+ n k \d_\Sp(Z, \rho(\g) Z).
\eqn
In particular 
\bqn
\d_\Sp(Z',\rho(\g) Z') \leq  n \d_\Sp(Z,\rho(\g)Z) + \frac{n-1}{k} \d_\Sp(Z',Z).
\eqn
Thus, for $A''= n$ and $B''>0$ fixed we can choose $k$ big enough such that 
$\frac{n-1}{k} \d_\Sp (Z',Z) \leq B''$ to get
\bqn
\d_\Sp(Z', \rho(\g) Z') \leq A'' \d_\Sp(Z, \rho(\g)Z) + B''.
\eqn

Since $\d_\Sp$ is left invariant, this implies 
\bqn
\d_\Sp(J(u), \rho(\g) J(u))&\leq&  \d_\Sp(\rho(\g)^N J(u), \rho(\g)^{N+1} J(u))  \\
&=& \d_\Sp(Z', \rho(\g) Z')\\
&\leq& A''\d_\Sp(Z, \rho(\g) Z)+B'', 
\eqn
hence the claim.
\end{proof}

\begin{lemma}\label{lem:compare}
There exist constants  $A'',B''>0$ depending only on $\dim(V)$ such that for all
$u\in T^1\DD$ and $\g \in \G_g$   
\bqn
\d_\Sp(J(u), \rho(\g) J(u)) \leq A'' \tr_\rho (\g) + 2B''.
\eqn
\end{lemma}
\begin{proof}
Fix some $\eps >0$. By Lemma~\ref{lem:subspace} there exists $Z_0\in \Yy_{L_-, s}$ such that 
\bqn
\d_\Sp(Z_0, \rho(\g) Z_0) \leq \inf_{X\in \Yy_{L_-, s}} \d_\Sp(X,
\rho(\g) X) +\eps = \tr_\rho(\g) +\eps. 
\eqn
Therefore Lemma~\ref{lem:lengths} implies 
\bqn
\d_\Sp (J(u),\rho(\g) J(u))\leq  A'' \d_\Sp(Z_0, \rho(\g)Z_0) + B''\leq A'' \tr_\rho(\g) +A''\eps+B''.
\eqn
Since this holds for all $\eps>0$ we get that 
\bqn
\d_\Sp (J(u), \rho(\g) J(u))\leq A'' \tr_\rho(\g) + 2B''.
\eqn
\end{proof}

\begin{proof}[Proof of Proposition~\ref{prop:symp}]
Let $\rho$ be a maximal representation of $\G_g$ into $\Sp(V) $ and let 
$p\in \DD$ be such that $ \tr_h(\g)=\d_\DD(p, \g p)$;,
then $p$ lies on the unique geodesic
$c_\g$ connecting $\g^-$ to $\g^+$.  

Let $u= (\g^-, u_0,\g^+) \in T^1\DD$ be the unit tangent vector to
$c_\g$ at $p$ and $J(u)\in \Yy_{L_-, L_+}$ the image of $u$ under
the mapping $ J :T^1\DD
\to \Xx_\Sp$. By Lemma~\ref{lem:strongquasi} there
exists a constant $A'$ such that 
\bqn
A'^{-1} \tr_h(\g) = A'^{-1} \d_\DD(p, \g p) \leq \d_\Sp(J(u), \rho(\g) J(u)). 
\eqn
Applying Lemma~\ref{lem:compare} to this, there exist constants $A'', B''>0$ such that 
\bqn
\tr_h(\g) \leq A'\d_\Sp(J(u), \rho(\g) J(u))\leq 
A''A' \tr_\rho(\g) + 2A'B''.
\eqn
This in combination with Lemma~\ref{lem:upperbound} finishes the proof.
\end{proof}

\vskip1cm

\bibliographystyle{amsplain}
\bibliography{refs}

\providecommand{\bysame}{\leavevmode\hbox to3em{\hrulefill}\thinspace}
\providecommand{\MR}{\relax\ifhmode\unskip\space\fi MR }
\providecommand{\MRhref}[2]{%
  \href{http://www.ams.org/mathscinet-getitem?mr=#1}{#2}
}
\providecommand{\href}[2]{#2}
\begin{thebibliography}{10}

\bibitem{Bradlow_GarciaPrada_Gothen}
S.~B. Bradlow, O.~Garc\'\i a~Prada, and P.~B. Gothen, \emph{Surface group
  representations in {${\rm PU}(p,q)$} and {H}iggs bundles}, J. Diff. Geom.
  \textbf{64} (2003), no.~1, 111--170.

\bibitem{Bradlow_GarciaPrada_Gothen_homotopy}
\bysame, \emph{Homotopy groups of moduli spaces of representations}, preprint,
  {\sf arXiv:math.AG/0506444}, June 2005.

\bibitem{Burger_Iozzi_Labourie_Wienhard}
M.~Burger, F.~Labourie A.~Iozzi, and A.~Wienhard, \emph{Maximal representations
  of surface groups: {S}ymplectic {A}nosov structures}, Pure and Applied
  Mathematics Quaterly. Special Issue: In Memory of Armand Borel. \textbf{1}
  (2005), no.~2, 555--601.

\bibitem{Burger_Iozzi_Wienhard_tol}
M.~Burger, A.~Iozzi, and A.~Wienhard, \emph{Surface group representations with
  maximal {T}oledo invariant}, preprint.

\bibitem{Burger_Iozzi_Wienhard_tight}
\bysame, \emph{Tight embeddings}, preprint.

\bibitem{Burger_Iozzi_Wienhard_ann}
\bysame, \emph{Surface group representations with maximal {T}oledo invariant},
  C. R. Acad. Sci. Paris, S\'er. I \textbf{336} (2003), 387--390.

\bibitem{Clerc_Orsted_2}
J.~L. Clerc and B.~{\O}rsted, \emph{The {G}romov norm of the {K}aehler class
  and the {M}aslov index}, Asian J. Math. \textbf{7} (2003), no.~2, 269--295.

\bibitem{Domic_Toledo}
A.~Domic and D.~Toledo, \emph{The {G}romov norm of the {K}\"ahler class of
  symmetric domains}, Math. Ann. \textbf{276} (1987), no.~3, 425--432.

\bibitem{Douady}
A.~Douady, \emph{L'espace de {T}eichm\"uller}, Travaux de Thurston sur les
  surfaces, Asterisque 66-67, Soci\'et\'e de math\'ematique de france, 1979,
  pp.~127--137.

\bibitem{Farb_Margalit}
B.~Farb and D.~Margalit, \emph{A primer on mapping class groups}, In
  preparation.

\bibitem{Goldman_survey}
W.~M. Goldman, \emph{Mapping class group dynamics on surface group
  representations}, Preprint 2005.

\bibitem{Goldman_thesis}
\bysame, \emph{Discontinuous groups and the {E}uler class}, Thesis, University
  of California at Berkeley, 1980.

\bibitem{Goldman_88}
W.~M. Goldman, \emph{Topological components of spaces of representations},
  Invent. Math. \textbf{93} (1988), no.~3, 557--607.

\bibitem{Gothen}
P.~B. Gothen, \emph{Components of spaces of representations and stable
  triples}, Topology \textbf{40} (2001), no.~4, 823--850.

\bibitem{Hernandez}
L.~{H}ern\`andez Lamoneda, \emph{Maximal representations of surface groups in
  bounded symmetric domains}, Trans. Amer. Math. Soc. \textbf{324} (1991),
  405--420.

\bibitem{Ivanov_mapping}
N.~V. Ivanov, \emph{Mapping class groups}, Handbook of geometric topology,
  North-Holland, Amsterdam, 2002, pp.~523--633.

\bibitem{Korevaar_Schoen}
N.~J. Korevaar and R.~M. Schoen, \emph{Sobolev spaces and harmonic maps for
  metric space targets}, Comm. Anal. Geom. \textbf{1} (1993), no.~3-4,
  561--659.

\bibitem{Labourie_anosov}
F.~Labourie, \emph{Anosov flows, surface groups and curves in projective
  space}, to appear Inventiones Mathematicae, {\sf arXiv:math.DG/0401230}.

\bibitem{Labourie_energy}
\bysame, \emph{Cross {R}atio, {A}nosov {R}epresentations and the {E}nergy
  {F}unctional on {T}eichm\"uller space}, Preprint 2005.

\bibitem{Labourie_crossratio}
\bysame, \emph{Crossratios, surface groups, ${SL}(n,\mathbb{R})$ and
  ${C}^{1,h}({S}^1)\rtimes{D}iff^h({S}^1)$}, Preprint, {\sf
  arXiv:math.DG/0502441}.

\bibitem{GarciaPrada_Gothen_Mundet}
P.~Gothen O.~Garc{\'{\i}}a-Prada and I.~Mundet i~Riera, \emph{Connected
  components of the representation variety for {${\rm Sp}(2n,{\mathbb R})$}},
  Preprint in preparation.

\bibitem{Satake_book}
I.~Satake, \emph{Algebraic structures of symmetric domains}, Kan\^o Memorial
  Lectures, vol.~4, Iwanami Shoten, Tokyo, 1980.

\bibitem{Toledo_89}
D.~Toledo, \emph{Representations of surface groups in complex hyperbolic
  space}, J. Diff. Geom. \textbf{29} (1989), no.~1, 125--133.

\bibitem{Wienhard_thesis}
A.~Wienhard, \emph{Bounded cohomology and geometry}, Ph.D. thesis,
  Universit\"at Bonn, 2004, Bonner Mathematische Schriften Nr. 368.

\bibitem{Wienhard_grenoble}
\bysame, \emph{A generalisation of {T}eichm\"uller space in the {H}ermitian
  context}, S\'eminaire de Th\'eorie Spectrale et G\'eom\'etrie Grenoble
  \textbf{22} (2004), 103--123.

\end{thebibliography}
\vskip1cm


\end{document}